\documentclass[leqno,11pt]{article}

\usepackage{amsthm, amsmath}
\usepackage{amsfonts}

\makeatletter
  
\def\@seccntformat#1{\csname the#1\endcsname.\quad}
  \makeatother

\makeatletter
  
   \@addtoreset{equation}{section}
\makeatother

\makeatletter

\renewcommand\section{\@startsection {section}{1}{\z@}%
            {-4.5ex \@plus -1ex \@minus -.2ex}%
           {1.5ex \@plus.2ex}%
            {\normalfont\large\bfseries}}

\renewcommand\subsection{\@startsection {subsection}{1}{\z@}%
            {-3.5ex \@plus -1ex \@minus -.2ex}%
            {1.5ex \@plus.2ex}%
            {\normalfont\normalsize\bfseries}}

\def\<{\left<} \def\>{\right>}

\newtheorem{theorem}{Theorem}[section]
\newtheorem{proposition}{Proposition}[section]
\newtheorem{corollary}{Corollary}[section]
\newtheorem{lemma}{Lemma}[section]
\newtheorem{remark}{Remark}[section]

\newtheorem{definition}[theorem]{Definition}
\def\proof{\noindent{\it Proof. }}
\def\bea{\begin{eqnarray} }
\def\eea{\end{eqnarray} }
\def\be{\begin{equation} }
\def\ee{\end{equation} }

\title{Classification of biharmonic $\mathcal{C}$-parallel Legendrian 
submanifolds in $7$-dimensional Sasakian space forms
\footnote{Tohoku Mathematical Journal {\bf 71} (2019), 157-169. A section is added at the end.}}
\author{Toru Sasahara}
\date{}

\begin{document}


\maketitle


\begin{abstract}
{\footnotesize 
In  [5],  D. Fetcu and C. Oniciuc
 presented the classification result for
biharmonic $\mathcal{C}$-parallel Legendrian submanifolds in $7$-dimensional Sasakian space forms.
However, it is incomplete. In this paper, all such submanifolds are explicitly determined.}
\end{abstract}





\section{Introduction}

In \cite[Theorem 5.1]{fe},  Fetcu and Oniciuc  presented the classification result for
proper biharmonic $\mathcal{C}$-parallel Legendrian submanifolds  
in $7$-dimensional Sasakian space forms.
The case $(2)$ of the theorem is proved by applying  Lemma 4.4 in \cite{ba}. 
However, the Lemma is wrong, and hence Fetcu and Oniciuc's classification is incomplete. 
This paper corrects errors in \cite{ba}, and moreover,  completes the classification.

Our main result is the following, which 
determines explicitly all proper biharmonic $\mathcal{C}$-parallel Legendrian submanifolds  
in $7$-dimensional Sasakian space forms.
\begin{theorem}\label{main}
Let $f: M^3\rightarrow N^{7}(\epsilon)$ be a $3$-dimensional
$\mathcal{C}$-parallel Legendrian submanifold in a $7$-dimensional Sasakian space form of 
constant $\varphi$-sectional curvature $\epsilon$. 
Then  $M^3$ is proper biharmonic  if and only if  either{\rm :}

{\rm (1)} $M^3$ is flat, 
$N^{7}(\epsilon)=S^{7}(\epsilon)$ with $\epsilon>-1/3$,
where $S^7(\epsilon)$ is a unit sphere in $\mathbb{C}^4$
equipped with its canonical and deformed Sasakian structures,  and $f(M^3)$ is an open part of 
\begin{align}
f(u, v, w)=&\biggl(\frac{\lambda}{\sqrt{\lambda^2+\alpha^{-1}}}\exp\Bigl(i\Bigl(\frac{1}{\alpha\lambda}u\Bigr)\Bigr), \nonumber\\
&\frac{1}{\sqrt{\alpha(c-a)(2c-a)}}\exp(-i(\lambda u-(c-a)v)),\nonumber\\
&\frac{1}{\sqrt{\alpha\rho_1(\rho_1+\rho_2)}}\exp(-i(\lambda u+cv+\rho_1 w)),\nonumber\\
&\frac{1}{\sqrt{\alpha\rho_2(\rho_1+\rho_2)}}\exp(-i(\lambda u+cv-\rho_2 w))\biggr),
\label{legen0}
\end{align} 
where $\alpha=4/(\epsilon+3)$, $\rho_{1, 2}=(\sqrt{4c(2c-a)+d^2}\pm d)/2$ and $\lambda$, $a$,
 $c$, $d$ are real constants given by
\be
\begin{cases}
(3\lambda^2-\alpha^{-1})(3\lambda^4-2(\epsilon+1)\lambda^2+\alpha^{-2})+
\lambda^4((a+c)^2+d^2)=0,\\
(a+c)(5\lambda^2+a^2+c^2-7\alpha^{-1}+4)+cd^2=0,\\
d(5\lambda^2+d^2+3c^2+ac-7\alpha^{-1}+4)=0,\\
\alpha^{-1}+\lambda^2+ac-c^2=0\nonumber
\end{cases}
\ee
such that $-1/\sqrt{\alpha}<\lambda<0$, $0<a\leq (\lambda^2-\alpha^{-1})/\lambda$, 
$a\geq d\geq 0$, $a>2c$, $\lambda^2\ne1/(3\alpha)${\rm ;} \ or

{\rm (2)} $M^3$ is non-flat, 
$N^{7}(\epsilon)=S^{7}(\epsilon)$ with
$\epsilon\geq(-7+8\sqrt{3})/13$
 and $f(M^3)$ is an open part of 
\be
f(x, {\bf y})=\Biggl(\sqrt{\frac{\mu^2}{\mu^2+1}}
e^{-\frac{i}{\mu}x}, \sqrt{\frac{1}{\mu^2+1}}e^{i\mu x}{\bf y}
\Biggr)
,\label{legen}
\ee  
where ${\bf y}=(y_1, y_2, y_3)$, $||{\bf y}||=1$ and 
\be\mu^2=
\begin{cases}
1 & (\epsilon=1)\\
\displaystyle\frac{4\epsilon+4\pm\sqrt{13\epsilon^2+14\epsilon-11}}{3(3+\epsilon)} & (\epsilon\ne 1).\label{mu}
\end{cases}
\ee
\end{theorem}


\begin{remark}
{\rm 
The flat case (1) of Theorem \ref{main} has been proved by Fetcu and Oniciuc in 
\cite[Theorem 5.1]{fe}.
However, they did not give
 the explicit representation of non-flat biharmonic $\mathcal{C}$-parallel Legendrian submanifolds
 in $S^7(\epsilon)$.}\end{remark}

\begin{remark}\label{remark1}
{\rm The immersion (\ref{legen0}) can be  rewritten as 
\be
f(u, v, w)=\left(z_1(u), z_2(u){\bf y}(v, w)\right),\nonumber
\ee
where $(z_1(u), z_2(u))$ is a Legendre curve with constant curvature $(\lambda^2-\alpha^{-1})/\lambda$ in $S^3(\epsilon)$ given by
$$(z_1(u), z_2(u))=\left(\frac{\lambda}{\sqrt{\lambda^2+\alpha^{-1}}}e^{i\frac{1}{\alpha\lambda}u}, \frac{1}{\sqrt{\alpha\lambda^2+1}}e^{-i\lambda u}\right)$$
and ${\bf y}(u, v)$ is a Legendrian surface in $S^5(\epsilon)$ given by
$${\bf y}(v, w)=
\left(\frac{\sqrt{\alpha\lambda^2+1}}{\sqrt{\alpha(c-a)(2c-a)}}e^{i(c-a)v},
\frac{\sqrt{\alpha\lambda^2+1}}{\sqrt{\alpha\rho_1(\rho_1+\rho_2)}}e^{-i(cv+\rho_1 w)},
\frac{\sqrt{\alpha\lambda^2+1}}{\sqrt{\alpha\rho_2(\rho_1+\rho_2)}}e^{-i(cv-\rho_2 w)}\right).$$
}
\end{remark}

\begin{remark}\label{remark2}
{\rm (i) For each fixed $x$,    (\ref{legen}) has constant Gauss curvature
$(\mu^2+1)/\alpha$ with respect to the induced metric from $S^7(\epsilon)$.
We can check that the surface is an integral $\mathcal{C}$-parallel surface in $S^7(\epsilon)$.

(ii) The curve $$z(x):=\Bigl(\sqrt{\dfrac{\mu^2}{\mu^2+1}}e^{-\frac{i}{\mu}x}, \sqrt{\dfrac{1}{\mu^2+1}}e^{i\mu x}\Bigr)$$ given in
 (\ref{legen}) is a Legendre curve with constant curvature $(\mu^2-1)/(\mu\sqrt{\alpha})$
 in $S^3(\epsilon)$. }
\end{remark}


\begin{remark}
{\rm  
(i) In \cite[Theorem 5.1]{fe}, it is stated that when $\epsilon=5/9$, $M^3$ is locally isometric to a product  $\gamma\times\bar{M}^2$, where $\gamma$ is a curve  of  constant curvature $1/\sqrt{2}$ in $S^7(5/9)$ and $\bar{M}^2$ is 
 a surface of constant Gauss curvature $4/3$.
However, $1/\sqrt{2}$ should be replaced by $2/3$ because 
$\gamma$ coincides with $z(x)$ 
in Remark \ref{remark2}.

(ii) The function $\lambda$ in the case (2) of  
\cite[Theorem 5.1]{fe} and the function $\mu$ in
 (\ref{mu}) are related by the equation $\mu^2=\alpha\lambda^2$.
 Hence, in view of Remark \ref{remark2}, the case $\epsilon=1$ and 
the case 
$\mu^2=(4\epsilon+4+\sqrt{13\epsilon^2+14\epsilon-11})/(3(3+\epsilon))$ with
$\epsilon>1$ in (2) of Theorem \ref{main} are missing from \cite[Theorem 5.1]{fe}. 
}
\end{remark}



Applying Theorem \ref{main}, we have the following result
which corrects \cite[Corollary 5.2]{fe}.
\begin{corollary}\label{cor}
Let $f: M^3\rightarrow S^{7}(1)$ be a 
$\mathcal{C}$-parallel Legendrian submanifold.
Then  $M^3$ is proper biharmonic if and only if  either{\rm :}

 {\rm (1)} $M^3$ is flat, and $f(M^3)$ is an open part of 
\begin{align}
f(u, v, w)=&\biggl(-\frac{1}{\sqrt{6}}\exp(-i\sqrt{5}u), \nonumber\\
&\frac{1}{\sqrt{6}}\exp\Bigl(i\Bigl(\dfrac{1}{\sqrt{5}}u-\dfrac{4\sqrt{3}}{\sqrt{10}}v\Bigr)\Bigr),\nonumber\\
&\frac{1}{\sqrt{6}}\exp\Bigl(i\Bigl(\dfrac{1}{\sqrt{5}}u+\dfrac{\sqrt{3}}{\sqrt{10}}v-\dfrac{3\sqrt{2}}{2}w\Bigr)\Bigr),\nonumber\\
&\frac{1}{\sqrt{2}}\exp\Bigl(i\Bigl(\dfrac{1}{\sqrt{5}}u+\dfrac{\sqrt{3}}{\sqrt{10}}v+\dfrac{\sqrt{2}}{2}w\Bigr)\Bigr)\biggr);\nonumber \ or
\end{align}

{\rm (2)} $M^3$ is non-flat, and $f(M^3)$ is an open part of 
\be f(x, {\bf y})=\dfrac{1}{\sqrt{2}}(e^{ix}, e^{-ix}{\bf y}),\label{example}\ee
where ${\bf y}=(y_1, y_2, y_3)$ and  $||{\bf y}||=1$. 
\end{corollary}
\begin{remark}
{\rm 
The flat case (1) of Corollary \ref{cor} has been proved in \cite[Corollary 5.2]{fe}.
However, the non-flat submanifold (\ref{example}) is missing from \cite[Corollary 5.2]{fe}.}
\end{remark}
\begin{remark}
{\rm The author classified proper biharmonic Legendrian surfaces in $5$-dimensional Sasakian space forms (see \cite{sa1} and \cite{sa}). Those surfaces are flat and $\mathcal{C}$-parallel.}
\end{remark}

In the last section, by the same argument as in the proof of Theorem \ref{main}, 
we determine explicitly all proper biharmonic parallel Lagrangian submanifolds in $3$-dimensional complex projective space.

\section{Preliminaries}
\subsection{Sasakian space forms}
A $(2n+1)$-dimensional    manifold $N^{2n+1}$ is called an {\it almost contact manifold}  if it
admits a unit vector field $\xi$, a one-form $\eta$ and a $(1, 1)$-tensor field $\varphi$ satisfying
\bea
\eta(\xi)=1, \quad \varphi^2=-I+\eta\otimes\xi.\nonumber
\eea
Every almost contact manifold admits a Riemannian metric $g$ satisfying
\bea
g(\varphi X, \varphi Y)=g(X, Y)-\eta(X)\eta(Y).\nonumber
\eea
The quadruplet $(\varphi, \xi, \eta, g)$ is called an {\it almost contact metric structure}.
An almost contact metric structure is said to be {\it normal}
if the tensor field $S$ defined by
\bea
S(X, Y)=\varphi^2[X, Y]+[\varphi X, \varphi Y]-\varphi[\varphi X, Y]-\varphi[X, \varphi Y]
+2d\eta(X, Y)\xi\nonumber
\eea
vanishes identically. A normal almost contact structure is said to be {\it Sasakian}
if it satisfies
\bea
d\eta(X, Y):=(1/2)\left(X(\eta(Y))-Y(\eta(X))-\eta([X, Y])\right)=g(X, \varphi Y).\nonumber
\eea

The tangent plane in $T_pN^{2n+1}$ which is invariant
under $\varphi$ is called a $\varphi$-$section$.
The sectional curvature of $\varphi$-section  is called the {\it $\varphi$-sectional curvature}.
Complete and connected  Sasakian manifolds of constant
$\varphi$-sectional curvature are called {\it Sasakian space forms}.
Denote Sasakian space forms of constant
$\varphi$-sectional curvature $\epsilon$ by $N^{2n+1}(\epsilon)$.

Let $S^{2n+1}\subset {\mathbb C}^{n+1}$ be the unit hypersphere centered at the origin.
Denote  by $z$ the position vector field of $S^{2n+1}$ in ${\mathbb C}^{n+1}$
and by $g_0$ the induced metric. Let $\xi_0=-Jz$, where $J$ is the usual complex structure of ${\mathbb C}^{n+1}$
which is   defined by $JX=iX$ for $X\in T{\mathbb C}^{n+1}$. Let $\eta_0$ be a  $1$-form defined by
$\eta_0(X)=g_0(\xi_0, X)$ and 
 $\varphi_0$ be the tensor field defined by $\varphi_0=s\circ J$, where $s: T_z{\mathbb C}^{n+1}
\rightarrow T_zS^{2n+1}$ denotes the orthogonal projection. 
Then, $(S^{2n+1}, \varphi_0, \xi_0, \eta_0, g_0)$ is a Sasakian space form of constant $\varphi$-sectional curvature
$1$.
If we put 
\be
\eta=\alpha\eta_0, \quad \xi=\alpha^{-1}\xi_0, \quad\varphi=\varphi_0, \quad g=\alpha g_0+\alpha(\alpha-1)\eta_0\otimes\eta_0 \nonumber
\ee
for a positive constant $\alpha$,  then $(S^{2n+1}, \varphi, \xi, \eta, g)$ is a Sasakian space form
of constant $\phi$ sectional curvature $\epsilon=(4/\alpha)-3>-3$. We denote it by $S^{2n+1}(\epsilon)$.
Tanno \cite{tan} showed that a simply connected Sasakian space form $N^{2n+1}(\epsilon)$ with $\epsilon>-3$
 is isomorphic to $S^{2n+1}(\epsilon)$; i.e., there exists a $C^{\infty}$-diffeomorphism which maps the structure tensors
 into the corresponding structure tensors.


\subsection{Legendrian submanifolds in Sasakian space forms}
Let $M^m$ be an $m$-dimensional submanifold $M$ in a Sasakian space form $N^{2n+1}(\epsilon)$. 
If  $\eta$ restricted to $M^m$ vanishes, then $M^m$ is called an $integral$ 
$submanifold$,
in particular if $m=n$, it is called  a $Legendrian$ $submanifold$.
In particular a Legendrian submanifold in a $3$-dimensional Sasakian space form is called
a {\it Legendre curve}. One can see that a curve $z(s)$ in $S^3(\epsilon)\subset\mathbb{C}^2$ is a
Legendre curve if and only if it satisfies ${\rm Re}(z^{\prime}(s), iz(s))=0$ identically in $\mathbb{C}^2$, where
 $(\cdot, \cdot)$
 is the standard Hermitian inner product on $\mathbb{C}^2$. 

We denote  the second fundamental
form, the shape operator and the normal
connection of a submanifold by $h$, $A$ and $D$, respectively.
The mean curvature vector field $H$ is defined by $H=(1/m){\rm Tr}\hskip3pt h$. 
If it vanishes identically, then $M^m$ is called a {\it minimal submanifold}.
In particular, if $h\equiv 0$,  then $M^m$ is called a {\it totally geodesic submanifold}.
A Legendrian submanifold in a Sasakian manifold  is parallel, i.e.,
satisfies $\bar\nabla h=0$ if and only if it is totally geodesic. Here, $\bar\nabla h$ is defined by
$$({\bar\nabla}_{X}h)(Y,Z)= D_X h(Y,Z) - h(\nabla_XY,Z) - h(Y,\nabla_X Z).$$
A Legendrian submanifold is called {\it $\mathcal{C}$-parallel} 
if $\bar\nabla h$ is parallel to $\xi$.

For a Legendrian submanifold $M$ in a Sasakian space form,   we have (cf. \cite{bl})
\be
A_{\xi}=0, \ \ \varphi h(X, Y)=-A_{\varphi Y}X, \ \ 
\<h(X, Y), \varphi Z\>=\<h(X, Z), \varphi Y\>\label{symmetric}
\ee
for any vector fields $X$, $Y$ and $Z$ tangent to $M$, where $\<\cdot, \cdot\>$ is the inner product.
%
We denote by $K_{ij}$ the sectional curvature determined by an orthonormal pair $\{X_i, X_j\}$. Then
from the equation of Gauss we have
\be
K_{ij}=(\epsilon+3)/4+\<h(X_i, X_i), h(X_j, X_j)\>-||h(X_i, X_j)||^2.\label{gauss}
\ee

The following Legendrian submanifolds 
can be regarded as the simplest Legendrian submanifolds next to totally geodesic ones in Sasakian space
forms. 
\begin{definition}\label{h-umbilical}
{\rm An $n$-dimensional  Legendrian submanifold $M^n$ in a  Sasakian space form is called
 {\it $H$-umbilical} if  
 every point has a neighborhood $V$ on which there 
 exists an  orthonormal frame field $\{e_1, \ldots, e_n \}$
 such that
  the  second fundamental form takes the following form:
\be
\begin{split}\label{h-um}
&h(e_1, e_1)=\lambda \varphi e_1, \quad h(e_2, e_2)=\cdots=h(e_n, e_n)=\mu \varphi e_1,\\
&h(e_1, e_j)=\mu \varphi e_j, \quad h(e_j, e_k)=0, \quad j \ne k, \quad j, k=2, \ldots, n,\nonumber
\end{split}
\ee 
where $\lambda$ and $\mu$ are some functions on $V$. }
\end{definition}

\begin{remark}{\rm
If in Definition \ref{h-umbilical} we assume
 that  the mean curvature vector field  is nowhere vanishing, then
$e_1=-\varphi H/||H||$ holds and hence it is a  globally defined differentiable vector field, and
$\lambda$ is also a globally defined differentiable function.
Moreover, 
 at each point $p$ of  $M^n$,  the shape operator $A_{JH}$ has only one eigenvalue 
 $\mu(p)$
 on $D(p)=\{X\in T_pM^n| \<X, JH\>=0\}$.  Since $\mu=(n||H||-\lambda)/(n-1)$ holds, 
 it is  also a globally defined differentiable function.}
 \end{remark}

\subsection{Biharmonic submanifolds}
Let $f:M^n \rightarrow N$ be a smooth map between  two Riemannian manifolds.
The {\it tension field}
$\tau(f)$ of $f$ is a section
of the vector bundle $f^{*}TN$
defined by
$$
\tau(f):=
\sum_{i=1}^{n}\left\{\nabla^{f}_{e_i}df(e_i)
-df(\nabla_{e_i}e_i)\right\},
$$
where $\nabla^f$, $\nabla$ and $\{e_i\}$ denote the induced connection, the connection of $M^n$ and
a local orthonormal basis of $M^n$, respectively.

A smooth map $f$ is called
a {\it harmonic map} if 
it is a critical 
point of the energy functional
$$
E(f)=\int_{\Omega}||df||^2dv
$$
over every compact domain $\Omega$ of $M^n$, where $dv$ is the volume form of $M^n$.
A smooth map $f$ is harmonic if and only if $\tau(f)=0$ at each point on $M^n$ (cf. \cite{es2}).


The bienergy functional $E_2(f)$ of $f$ over compact domain $\Omega\subset M^n$ is defined by 
\be
E_2(f)=\int_{\Omega}||\tau(f)||^2dv.\nonumber
\ee
Thus $E_2$ provides a measure for the extent to which $f$ fails to be harmonic. 
If $f$ is a critical point of $E_2$
over every compact domain $\Omega$, then $f$ is called a
{\it biharmonic map}. 
In \cite{ji2}, Jiang proved that $f$ is biharmonic if and only if its bitension field defined by
\be
\tau_2(f):=\sum_{i=1}^{n}\left\{(
\nabla^{f}_{e_i}\nabla^{f}_{e_i}-\nabla^{f}_
{\nabla_{e_i}e_i})\tau(f)+R^{N}(\tau(f),df(e_i))df(e_i)\right\}\nonumber
\ee 
vanishes identically, 
where $R^{N}$ is 
the curvature tensor  of $N$. 

A submanifold is called a {\it biharmonic submanifold} if the isometric immersion that defines the submanifold
is biharmonic map. Minimal submanifolds are biharmonic.
A biharmonic submanifold is said to be a {\it proper} biharmonic submanifold if it is non-minimal.

Loubeau and Montaldo introduced a  class which includes biharmonic submanifolds as follows.
\begin{definition}[\cite{LM3}]
{\rm
An isometric immersion $f: M\to N$
is called {\it biminimal\/} 
if it is a critical point of
the bienergy functional $E_2$ 
with respect to all
{\it normal variation} with compact support.
Here, a normal variation means a 
variation $f_t$ through $f=f_0$ 
such that the variational vector field 
$V=df_{t}/dt|_{t=0}$ is normal to $f(M)$. 
In this case, $M$ or $f(M)$ is called a {\it biminimal submanifold} in $N$.}
\end{definition}

An isometric  immersion $f$ is biminimal  if and only if the normal part of $\tau_2(f)$ vanishes identically.
Clearly, biharmonic submanifolds are biminimal. 
Biminimal $H$-umbilical 
Legendrian submanifolds in Sasakian space forms have been classified by the author as follows.

\begin{theorem}[\cite{sa}]\label{biminimal}
Let $f: M^n\rightarrow N^{2n+1}(\epsilon)$ be a 
non-minimal biminimal $H$-umbilical Legendrian submanifold, where $n\geq 3$.
 Then $N^{2n+1}(\epsilon)=S^{2n+1}(\epsilon)$ with
\be
\epsilon\geq\frac{- 3 n^2-2n+5+32\sqrt{n}}{n^2+6n+25} \hskip 5pt (>-3) \nonumber
\ee and $f(M^n)$ is an open part of 
\be
f(x, {\bf y})=\Biggl(\sqrt{\frac{\mu^2}{\mu^2+1}}
e^{-\frac{i}{\mu}x}, \sqrt{\frac{1}{\mu^2+1}}e^{i\mu x}{\bf y}
\Biggr),\nonumber
\ee  
where ${\bf y}=(y_1, \ldots, y_{n})$, $||{\bf y}||=1$ and 
\be\mu^2=
\begin{cases}
1 & (\epsilon=1)\\
\displaystyle\frac{(n+5)\epsilon+3n-1\pm\sqrt{P(n, \epsilon)}}{2(3+\epsilon)n} & (\epsilon\ne 1),\nonumber
\end{cases}
\ee
where $P(n, \epsilon):=(n^2+6n+25)\epsilon^2+(6n^2+4n-10)\epsilon+9 n^2-42n+1$.
\end{theorem}

\begin{remark}
{\rm Submanifolds given in Theorem \ref{biminimal} are in fact proper biharmonic.}
\end{remark}

\section{$\mathcal{C}$-parallel Legendrian submanifolds}
\subsection{A special orthonormal basis}
We recall a special local orthonormal basis which is  used in \cite{ba} (see also \cite{fe}).
Let $M$ be a non-minimal Legendrian submanifold of $N^7(\epsilon)$. Let $p$ be an 
arbitrary point of $M$, and denote by $U_pM$  the unit sphere in $T_pM$.  We consider the function 
$f_p: U_pM\rightarrow \mathbb{R}$ given by
\be
f_p(u)=\<h(u , u), \varphi u\>.\nonumber
\ee  
A function $f_p$ attains a critical value at $X$ if and only if $\<h(X, X),
 \varphi Y\>=0$ for all $Y\in U_pM$ with $\<X, Y\>=0$, i.e., $X$ is an eigenvector of $A_{\varphi X}$.

We take $X_1$ as a vector at which $f_p$ attains its maximum. 
Then there exists a local orthonormal basis $\{X_1, X_2, X_3\}$ of $T_pM$ such that the shape operators take the following forms (cf. \cite{ba}): 
\be
A_{\varphi X_1}= \left(
    \begin{array}{ccc}
      \lambda_1 & 0 & 0 \\
      0 & \lambda_2 & 0\\ 
      0 & 0 & \lambda_3 \\
   \end{array}
  \right),\ \  
  A_{\varphi X_2}=\left(
    \begin{array}{ccc}
      0& \lambda_2 & 0 \\
      \lambda_2 & a & b\\ 
      0 & b & c \\
   \end{array}
  \right),\ \ 
  A_{\varphi X_3}=\left(
    \begin{array}{ccc}
      0 & 0 & \lambda_3 \\
      0 & b & c\\ 
      \lambda_3 & c & d \\
   \end{array}
  \right),\label{shape}
\ee
where
\be
\lambda_1>0,  \ \ \lambda_1\geq 2\lambda_2, \ \ 
\lambda_1\geq 2\lambda_3,\ \ a\geq 0, \ \ a^2\geq d^2,\label{lambda}
\ee
and moreover, if $\lambda_2=\lambda_3$, then $b=0$ and $a\geq 2c$.



\begin{lemma}\label{special}
The vector  $X_1\in T_pM$ can be differentiably extended to a vector field $X_1(x)$ 
on a neighborhood $V$ of $p$ such that at
every point $x$ of $V$, $f_x$ attains a critical value at $X_1(x)$, that is, 
$X_1(x)$ is an eigenvector of $A_{\varphi X_1(x)}$.
\end{lemma}
\proof
Let $E_1(x)$, $E_2(x)$, $E_3(x)$ be an arbitrary local 
differentiable orthonormal frame field on a neighborhood $V$ of $p$, such that $E_i(p)=X_i$.
The purpose is to find a local differentiable vector field $X_1(x)=\sum y^i(x)E_i(x)$ such that $(y^1)^2+(y^2)^2+(y^3)^2=1$ and at every point $x$ of $V$, $f_x$
 attains a critical value at $X_1(x)$. 
As in the proof of Theorem A in \cite{li}, we apply Lagrange's multiplier method.

Consider the following function:
\be
F(x, y^1, y^2, y^3, \lambda):=\sum_{i, j, k} h_{ijk}y^iy^jy^k-\lambda\{(y^1)^2+(y^2)^2+(y^3)^2-1\},\nonumber
\ee
where $h_{ijk}:=\<h(E_i(x), E_j(x)), \varphi E_k(x)\>$.
We shall show that there exist differentiable functions $y^1$, $y^2$, $y^3$ defined a neighborhood of $p$ satisfying 
the following system of equations: 
\be
\begin{cases}
\displaystyle\frac{\partial F}{\partial y^i}=3\sum_{j, k}h_{ijk}(x)y^jy^k-2y^i\lambda=0, \ \ i\in\{1, 2, 3\},\\
\displaystyle\frac{\partial F}{\partial\lambda}=(y^1)^2+(y^2)^2+(y^3)^2-1=0.\label{system}
\end{cases}
\ee
Define functions $G_i$ by
\be
\begin{cases}
G_i(x, y^1, y^2, y^3, \lambda)=3\displaystyle\sum_{j, k} h_{ijk}(x)y^jy^k-2y^i\lambda \ \ {\rm for}\ \ i=1, 2, 3.\nonumber\\
G_{4}(x, y^1, y^2, y^3)=(y^1)^2+(y^2)^2+(y^3)^2-1.
\end{cases}
\ee
Since $X_1(p)=X_1=E_1(p)$, we have $(y^1, y^2, y^3)=(1, 0, 0)$ at $p$.
It follows from (\ref{shape}) and (\ref{system}) that
$2\lambda(p)=3\lambda_1$.  We set $y^{4}=\lambda$. A straightforward computation yields
\be
\det\biggl(\frac{\partial G_\alpha}{\partial y^\beta}\biggr)(p)=36(\lambda_2-\lambda_1)(\lambda_3-\lambda_1).\label{det}
\ee
By (\ref{lambda}), we have 
$\lambda_2\ne\lambda_1\ne\lambda_3$. Hence the RHS of (\ref{det}) is not zero. The implicit function theorem shows that there exist local differentiable functions $y^1(x)$, $y^2(x)$, $y^3(x)$, $\lambda(x)$ on a neighborhood of $p$ satisfying (\ref{system}). The proof is finished.
\qed

\begin{remark}
{\rm In \cite{ba} and \cite{fe}, the differentiablity of $X_1(x)$ is not proved.}
\end{remark}

If the eigenvalues of $A_{\varphi X_1(x)}$ have constant multiplicities on a neighborhood $V$ of $p$, we can extend $X_2$ and $X_3$ differentiably to vector fields $X_2(x)$ and $X_3(x)$
on $V$. We work on the open dense set of $M$ defined by this property.

\subsection{Correction to a paper by Biakoussis, Blair and Koufogiorgos}
 Let $M$ be a $\mathcal{C}$-parallel Legendrian  submanifold of $N^7(\epsilon)$. 
The condition that $M$ is $\mathcal{C}$-parallel is equivalent to  $\nabla\varphi h=0$, where $\nabla$
is the Levi-Civita connection of $M$. 
Hence we have 
\be
R\cdot\varphi h=0,\label{semi}
\ee
where $R$ is the curvature tensor of $M$.

By using (\ref{gauss}), (\ref{shape}) and (\ref{semi}),  
Biakoussis et al. obtained
a system of algebraic equations  with respect to $\lambda_1$, $\lambda_2$, $\lambda_3$, $a$, $b$, $c$, $d$, $K_{12}$, $K_{13}$ and $K_{23}$ (see \cite[pages 211-212]{ba}).

However, the equation $(3.19)$-(iv) in \cite{ba}, i.e., $c(a-2c)(\lambda_2-\lambda_3)=0$ is incorrect. 
It should be replaced by $$b(a-2c)(\lambda_2-\lambda_3)=0,$$
which is obtained by $\<(R(X_1, X_2)\cdot\varphi h)(X_2, X_2), X_3\>=0.$

In \cite[Lemma 4.4]{ba}, it is
stated that if  $\lambda_1=2\lambda_3\ne 2\lambda_2$, then  $\epsilon=-3$.
However, the proof is based on the the wrong equation
$(3.19)$-(iv) (see page 214, line 11), and hence the statement is also wrong.
The following is a counterexample to \cite[Lemma 4.4]{ba}:
The submanifold (\ref{example}) is a $H$-umbilical Legendrian submanifold
such that, with respect to some orthonormal local frame field $e_1, e_2, e_3$ with
$e_1=\partial/\partial x$, the second fundamental form $h$ satisfies
\be
\begin{split}
&h(e_1, e_1)=0, \ \  h(e_2, e_2)=h(e_3, e_3)=\varphi e_1,\\
&h(e_1, e_2)=\varphi e_2, \ \ h(e_1, e_3)=\varphi e_3, \ \  h(e_2, e_3)=0.\nonumber
\end{split}
\ee 
We put
 $X_1=(e_1-\sqrt{2}e_2)/\sqrt{3}$, $X_2=(\sqrt{2}e_1+e_2)/\sqrt{3}$ and $X_3=e_3$.
Then the shape operators  take the forms (\ref{shape}) with $\lambda_1=2/\sqrt{3}$, $\lambda_2=-1/\sqrt{3}$, $\lambda_3=1/\sqrt{3}$, $a=c=\sqrt{2}/\sqrt{3}$ and $b=d=0$.

On the other hand, 
following the wrong statement of \cite[Lemma 4.4]{ba}, the non-flat case (2)  of \cite[Theorem 5.1]{fe}  is investigated. 
Therefore, the classification presented in the theorem is incomplete.


\subsection{Biharmonic $\mathcal{C}$-parallel Legendrian submanifolds}

We shall prove Theorem \ref{main}.
First, we recall the following.
\begin{proposition}[\cite{fe}]\label{epsilon}
A $\mathcal{C}$-parallel Legendrian submanifolds in  a $7$-dimensional Sasakian space form $N^7(\epsilon)$ is proper biharmonic if and only if  $\epsilon>-1/3$ and 
\be
\mathrm{Tr}\hskip2pt h(\cdot, A_H\cdot)=(3\epsilon+1)/2.\label{biharmonic}
\ee
\end{proposition}
By applying the proof of   \cite[Lemmas 4.2-4.6]{ba} and Proposition \ref{epsilon}, we obtain the following.

\begin{proposition}\label{umbilical}
Let $M^3$ be a proper biharmonic   $\mathcal{C}$-parallel Legendrian submanifold in $N^7(\epsilon)$.
If $M$ is non-flat, then it is $H$-umbilical.
\end{proposition}
\proof
By \cite[Lemma 4.2]{ba}, the case $\lambda_1\ne2\lambda_2\ne 2\lambda_3\ne\lambda_1$ cannot hold.
According to the proof of  \cite[Lemma 4.6]{ba}, 
the case $\lambda_1=2\lambda_2=2\lambda_3$ cannot hold for $\epsilon>-3$. 
Hence, by Proposition \ref{epsilon} the proof is  divided into the following three cases.

{\bf Case (i).} $\lambda_1=2\lambda_2\ne 2\lambda_3$.
In the proof of \cite[Lemma 4.3]{ba}, we have 
\be
\lambda_1=2\lambda_2=-\lambda_3=\sqrt{2(\epsilon+3)}/4, \ \ a=c=d=0, \ \ b=\pm\sqrt{6(\epsilon+3)}/8. \label{case1.1}
\ee
We choose a local orthonormal frame field $\{e_1, e_2 ,e_3\}$ as follows:
\be
e_1=(X_1\pm\sqrt{3}X_3)/2, \ \ e_2=X_2, \ \ e_3=(\mp\sqrt{3}X_1+X_3)/2,\nonumber
\ee
where the $\pm$ signs are determined by the sign of $b$.
Then, by a straightforward computation using (\ref{case1.1}), we obtain
\be
\begin{split}
&h(e_1, e_1)=-(\sqrt{2(\epsilon+3)}/4)\varphi e_1, \ \  h(e_2, e_2)=h(e_3, e_3)
=(\sqrt{2(\epsilon+3)}/4)\varphi e_1,\\
&h(e_2, e_3)=0, \ \ h(e_1, e_i)=(\sqrt{2(\epsilon+3)}/4)\varphi e_i,  \ \ i\in\{2, 3\},\label{case1}
\end{split}
\ee
which implies that $M$ is $H$-umbilical. Moreover, from (\ref{biharmonic}) and (\ref{case1}) we have
$\epsilon=5/9$ (see the subcase (a) of (2) in \cite[Theorem 5.1]{fe}).

{\bf Case (ii).} $\lambda_1=2\lambda_3\ne 2\lambda_2$.
Following the proof of \cite[Lemma 4.4]{ba} (page 214, lines 7-10),
we have
\begin{align}
& K_{12}=0, \label{K3}\\ 
& b=d=0, \ \ c\ne 0. \label{bcd}
\end{align}
Moreover, in \cite[(3.16)-(iv), (3.21)]{ba}) the following equations have been obtained:
\begin{align}
& c(K_{23}+\lambda_3(\lambda_2-\lambda_3))=0, \label{K1}\\
 & (\lambda_2-\lambda_3)(K_{23}-b^2-c^2)=0.\label{K2}
\end{align}
From  (\ref{bcd}), (\ref{K1}), (\ref{K2}) and $\lambda_2\ne\lambda_3$, we have
\be
 K_{23}+\lambda_3(\lambda_2-\lambda_3)=0, \ \ K_{23}=c^2.\label{K4}
\ee
We note that
 (\ref{bcd}) and (\ref{K4}) show
$\lambda_3\ne0$.
It follows from (\ref{gauss}), (\ref{shape}), (\ref{K3}) and (\ref{K4}) that
\be
\begin{cases}\label{lam}
\lambda_2^2=4c^2-2ac-\beta,\\
\lambda_3^2=3c^2-ac-\beta,\\
\lambda_2\lambda_3=2c^2-ac-\beta,
\end{cases}
\ee
where $\beta=(\epsilon+3)/4$. 

We choose a local orthonormal frame field $\{e_1, e_2 ,e_3\}$ as follows:
\be
e_1=(\lambda_3X_1+cX_2)/\sqrt{\lambda_3^2+c^2}, \ \
e_2=(-cX_1+\lambda_3X_2)/\sqrt{\lambda_3^2+c^2}, \ \ e_3=X_3.\nonumber
\ee
We set
\be k(a, c):=8c^4-6ac^3+(a^2-3\beta)c^2+a\beta c.\nonumber
\ee 
Then, by a straightforward computation using (\ref{lam}), we obtain  
\be
\begin{split}\label{h}
&\<h(e_1, e_1), \varphi e_2\>=
\frac{ck(a, c)}{\lambda_3(\lambda_3^2+c^2)^{3/2}}, \ \ \<h(e_1, e_1), \varphi e_3\>=0, \\
 & \<h(e_2, e_2)-h(e_3, e_3), \varphi e_1\>=-\frac{k(a, c)}{(\lambda_3^2+c^2)^{3/2}},\ \ 
h(e_2, e_3)=0, \\
& \<h(e_2, e_2), \varphi e_2\>=-\frac{\lambda_3 k(a,c)}{c(\lambda_3^2+c^2)^{3/2}},\ \ \<h(e_3, e_3), \varphi e_3\>=0.
\end{split}
\ee

On the other hand,  substituting (\ref{lam}) into the identity  $\lambda_2^2\lambda_3^3- (\lambda_2\lambda_3)^2=0$ gives
\be
k(a, c)=0.\nonumber
 \ee
Hence, it follows from (\ref{symmetric}) and (\ref{h}) that $M$ is $H$-umbilical.

{\bf Case (iii).} $\lambda_1\ne 2\lambda_2=2\lambda_3$.
By rotating the vector fields $X_2$ and $X_3$, if necessary, we may assume that $b=0$.
In \cite[Lemma 4.5]{ba}, it is  proved that 
if $M$ is non-flat,  then $a\ne 2c$ and $a=c=d=0$. Thus,  $M$ is $H$-umbilical.
\qed

\medskip
{\bf Proof of Theorem \ref{main}}: 
The flat case (1) has been proved in (1) of  \cite[Theorem 5.1]{fe}.
 Applying Proposition \ref{umbilical} and
 Theorem \ref{biminimal} for $n=3$,  we can prove the non-flat case 
 (2).\qed

\begin{remark}
{\rm In \cite{fe}, the case (ii) of Proposition \ref{umbilical} was not investigated.}
\end{remark}

\section{Biharmonic parallel Lagrangian submanifolds}

Let  $\mathbb{C}P^n(4)$ denote the complex projective space  of complex dimension $n$ and constant
holomorphic  sectional curvature $4$.  We denote by $J$ the almost 
complex structure of $\mathbb{C}P^n(4)$. 
An $n$-dimensional submanifold $M^n$ of $\mathbb{C}P^n(4)$ is said to be {\it Lagrangian} if $J$ interchanges the tangent and the normal spaces at each point.

In \cite[Theorem 6.3]{fe}, Fetcu and Oniciuc presented the classification result of 
 proper biharmonic parallel Lagrangian submanifolds in  
$\mathbb{C}P^3(4)$. However, the theorem 
 is proved by applying the wrong statement of \cite[Lemma 4.4]{ba}, and hence the classification
 is incomplete.
This section completes it.
 First, we recall the following.
\begin{proposition}[\cite{fe}]\label{parallel}
Let $L: M^3\rightarrow \mathbb{C}P^3(4)$ be a proper biharmonic parallel Lagrangian immersion. 
Then $L$ is locally given by $\pi\circ f$, where $\pi: S^{2n+1}(1)\rightarrow \mathbb{C}P^n(4)$ is the Hopf fibration and $f: M^3\rightarrow S^7(1)$ is a non-minimal $\mathcal{C}$-parallel Legendrian immersion satisfying
 \be
 \mathrm{Tr}\hskip2pt h(\cdot, A_H\cdot)=6H.\label{-4}\nonumber
 \ee   
\end{proposition}

The following theorem determines explicitly all
 proper biharmonic parallel Lagrangian submanifolds in  
$\mathbb{C}P^3(4)$.
\begin{theorem}\label{theorem2}
Let $L: M^3\rightarrow \mathbb{C}P^3(4)$ a proper biharmonic parallel Lagrangian submanifold.
Then $L$ is locally congruent to $\pi\circ f$, where $f: M^3\rightarrow S^{7}(1)$ is 
one of the following{\rm :}

{\rm (1)} $M^3$ is flat and
\begin{align}
f(u, v, w)=&\biggl(\frac{\lambda}{\sqrt{\lambda^2+1}}\exp\Bigl(i\Bigl(\frac{1}{\lambda}u\Bigr)\Bigr), \nonumber\\
&\frac{1}{\sqrt{(c-a)(2c-a)}}\exp(-i(\lambda u-(c-a)v)),\nonumber\\
&\frac{1}{\sqrt{\rho_1(\rho_1+\rho_2)}}\exp(-i(\lambda u+cv+\rho_1 w)),\nonumber\\
&\frac{1}{\sqrt{\rho_2(\rho_1+\rho_2)}}\exp(-i(\lambda u+cv-\rho_2 w))\biggr),\label{legen3}
\end{align} 
where $\rho_{1, 2}=(\sqrt{4c(2c-a)+d^2}\pm d)/2$ and 
the $4$-tuple $(\lambda, a, c, d)$ is given by one of the following{\rm :}
\begin{align*}
&\left(-\sqrt{\frac{4-\sqrt{13}}{3}}, \sqrt{\frac{7-\sqrt{13}}{6}}, -\sqrt{\frac{7-\sqrt{13}}{6}}, 0\right),\\
&\left(-\sqrt{\frac{1}{5+2\sqrt{3}}}, \sqrt{\frac{45+21\sqrt{3}}{13}},
 -\sqrt{\frac{6}{21+11\sqrt{3}}}, 0\right),\\
 &\left(-\sqrt{\frac{1}{6+\sqrt{13}}}, \sqrt{\frac{523+139\sqrt{13}}{138}}, -\sqrt{\frac{79-17\sqrt{13}}{138}},
 \sqrt{\frac{14+2\sqrt{13}}{3}}\right){\rm ;}
\end{align*}

{\rm (2)} $M^3$ is non-flat and
\be
f(x, {\bf y})=\Biggl(\sqrt{\frac{\mu^2}{\mu^2+1}}
e^{-\frac{i}{\mu}x}, \sqrt{\frac{1}{\mu^2+1}}e^{i\mu x}{\bf y}
\Biggr),\label{legen2}
\ee  
where ${\bf y}=(y_1, y_2, y_3)$, $||{\bf y}||=1$ and 
$\mu^2=(4\pm\sqrt{13})/3$.
\end{theorem}
\proof
The flat case (1) has been proved in \cite[Corollary 6.4]{fe}.
Applying Propositions \ref{umbilical} and  \ref{parallel} and modifying
the second equation of \cite[(5.33)]{sa} to $\lambda^2+2\mu^2=6$,
we can prove the non-flat case (2).\qed

\begin{remark}
{\rm Fetcu and Oniciuc \cite{fe} did not give
 the explicit representation of non-flat proper 
 biharmonic parallel Lagrangian submanifolds in $\mathbb{C}P^3(4)$.}
\end{remark}


\begin{remark}
{\rm The immersion (\ref{legen3})  can be  rewritten as
 the one with $\alpha=1$ in  Remark 
\ref{remark1} (cf. 
\cite{dillen}, \cite{li2}).
 }
\end{remark}

\begin{remark}
{\rm The  immersion (\ref{legen2}) has 
the same properties as in Remark \ref{remark2}, where $\alpha=1$. 
From this, we see that  (\ref{legen2}) with  $\mu^2=(4+\sqrt{13})/3$ is missing from \cite[Theorem 6.3]{fe}. 
}
\end{remark}

\begin{remark}
{\rm The author classified proper biharmonic Lagrangian surfaces of constant mean curvature
 in $\mathbb{C}P^2(4)$ (see \cite{sa2}). Those surfaces are flat and parallel.}
\end{remark}

\section{Corrections to this paper {\normalsize\rm{(added on November 16, 2022)}}}
\subsection{Correction to Lemma \ref{special}}
Equation (\ref{det}) should be replaced by
\be
\det\biggl(\frac{\partial G_\alpha}{\partial y^\beta}\biggr)(p)=36(2\lambda_2-\lambda_1)(2\lambda_3-\lambda_1).\label{det2}
\ee
Therefore, Lemma \ref{special} should be replaced by

\begin{lemma}\label{lemma2}
If $\lambda_1\ne 2\lambda_2$ and $\lambda_1\ne 2\lambda_3$, then the vector  $X_1\in T_pM$ can be differentiably extended to a vector field $X_1(x)$ 
on a neighborhood $V$ of $p$ such that at
every point $x$ of $V$, $f_x$ attains a critical value at $X_1(x)$, that is, 
$X_1(x)$ is an eigenvector of $A_{\varphi X_1(x)}$.
\end{lemma}

\subsection{Correction to the proof of Proposition \ref{umbilical}} 
Proof of  Proposition \ref{umbilical}: 
Let $p$ be an arbitrary point of $M$, and we choose a local orthonormal basis
$\{X_1, X_2, X_3\}$ of $T_pM$ such that the shape operators take the form (\ref{shape}) with (\ref{lambda}).

By an argument given in the proof of \cite[Lemma 4.2]{ba}, the case $\lambda_1\ne2\lambda_2\ne 2\lambda_3\ne\lambda_1$ cannot hold.
According to the proof of  \cite[Lemma 4.6]{ba}, 
the case $\lambda_1=2\lambda_2=2\lambda_3$ cannot hold for $\epsilon>-3$. 
Note that these two assertions can be obtained without using (3.4)-(3.13) in \cite{ba}.
By Proposition \ref{epsilon}, 
$\epsilon>-1/3$ must be satisfied, and hence the proof is  divided into the following three cases.

{\bf Case (i).} $\lambda_1=2\lambda_2\ne 2\lambda_3$.
In the proof of \cite[Lemma 4.3]{ba}, we have 
\be
\lambda_1=2\lambda_2=-\lambda_3=\sqrt{2(\epsilon+3)}/4, \ \ a=c=d=0, \ \ b=\pm\sqrt{6(\epsilon+3)}/8. \label{case1.11}
\ee
We choose a local orthonormal basis $\{e_1, e_2 ,e_3\}$ of $T_pM$ as follows:
\be
e_1=(X_1\pm\sqrt{3}X_3)/2, \ \ e_2=X_2, \ \ e_3=(\mp\sqrt{3}X_1+X_3)/2,\nonumber
\ee
where the $\pm$ signs are determined by the sign of $b$.
Then, by a straightforward computation using (\ref{case1.11}), we obtain
\be
\begin{split}
&h(e_1, e_1)=-(\sqrt{2(\epsilon+3)}/4)\varphi e_1, \ \  h(e_2, e_2)=h(e_3, e_3)
=(\sqrt{2(\epsilon+3)}/4)\varphi e_1,\\
&h(e_2, e_3)=0, \ \ h(e_1, e_i)=(\sqrt{2(\epsilon+3)}/4)\varphi e_i,  \ \ i\in\{2, 3\}.\label{case11}
\end{split}
\ee

{\bf Case (ii).} $\lambda_1=2\lambda_3\ne 2\lambda_2$.
Following the proof of \cite[Lemma 4.4]{ba} (page 214, lines 7-10),
we have
\begin{align}
& K_{12}=0, \label{K33}\\ 
& b=d=0, \ \ c\ne 0. \label{bcdd}
\end{align}
Moreover, in \cite[(3.16)-(iv), (3.21)]{ba}) the following equations have been obtained:
\begin{align}
& c(K_{23}+\lambda_3(\lambda_2-\lambda_3))=0, \label{K11}\\
 & (\lambda_2-\lambda_3)(K_{23}-b^2-c^2)=0.\label{K22}
\end{align}
From  (\ref{bcdd}), (\ref{K11}), (\ref{K22}) and $\lambda_2\ne\lambda_3$, we have
\be
 K_{23}+\lambda_3(\lambda_2-\lambda_3)=0, \ \ K_{23}=c^2.\label{K44}
\ee
We note that
 (\ref{bcdd}) and (\ref{K44}) show
$\lambda_3\ne0$.
It follows from (\ref{gauss}), (\ref{shape}), (\ref{K33}) and (\ref{K44}) that
\be
\begin{cases}\label{lamm}
\lambda_2^2=4c^2-2ac-\beta,\\
\lambda_3^2=3c^2-ac-\beta,\\
\lambda_2\lambda_3=2c^2-ac-\beta,
\end{cases}
\ee
where $\beta=(\epsilon+3)/4$. 

We choose a local orthonormal basis $\{e_1, e_2 ,e_3\}$ of $T_pM$ as follows:
\be
e_1=(\lambda_3X_1+cX_2)/\sqrt{\lambda_3^2+c^2}, \ \
e_2=(-cX_1+\lambda_3X_2)/\sqrt{\lambda_3^2+c^2}, \ \ e_3=X_3.\nonumber
\ee
We set
\be k(a, c):=8c^4-6ac^3+(a^2-3\beta)c^2+a\beta c.\nonumber
\ee 
Then, by a straightforward computation using (\ref{lamm}), we obtain  
\be
\begin{split}\label{hh}
&\<h(e_1, e_1), \varphi e_2\>=
\frac{ck(a, c)}{\lambda_3(\lambda_3^2+c^2)^{3/2}}, \ \ \<h(e_1, e_1), \varphi e_3\>=0, \\
 & \<h(e_2, e_2)-h(e_3, e_3), \varphi e_1\>=-\frac{k(a, c)}{(\lambda_3^2+c^2)^{3/2}},\ \ 
h(e_2, e_3)=0, \\
& \<h(e_2, e_2), \varphi e_2\>=-\frac{\lambda_3 k(a,c)}{c(\lambda_3^2+c^2)^{3/2}},\ \ \<h(e_3, e_3), \varphi e_3\>=0.
\end{split}
\ee

On the other hand,  substituting (\ref{lamm}) into the identity  $\lambda_2^2\lambda_3^3- (\lambda_2\lambda_3)^2=0$ gives
\be
k(a, c)=0.\nonumber
 \ee
Hence, it follows from (\ref{symmetric}) and (\ref{hh}) that the second fundamental form takes 
the form in Definition \ref{h-umbilical} at $p$.

{\bf Case (iii).} $\lambda_1\ne 2\lambda_2=2\lambda_3$. By Lemma \ref{lemma2}, the basis $\{X_1, X_2, X_3\}$ of $T_pM$
can be differentiably extended to an orthonormal frame field $\{X_1(x), X_2(x), X_3(x)\}$ 
on a neighborhood $V$ of $p$ such that at
every point $x$ of $V$ the shape operators take the form (\ref{shape}) with $\lambda_1\ne 2\lambda_2$ and 
$\lambda_1\ne 2\lambda_3$.
It follows from  (3.4)-(3.10) in \cite{ba} that  $\lambda_1$, $\lambda_2$, $\lambda_3$ and $a$ are constant on $V$.
By rotating the vector fields $X_2(x)$ and $X_3(x)$, if necessary, we may assume that $b=0$ on $V$.
In \cite[Lemma 4.5]{ba}, it is proved that 
if $M$ is non-flat,  then $a=c=d=0$ on $V$. 

Consequently, $M$ is $H$-umbilical.
\qed

 \vskip20pt
 
 
 Center for Liberal Arts and Sciences
 
Hachinohe Institute of Technology
 
Hachinohe 031-8501

JAPAN

{{\it E-mail address}: {\tt sasahara@hi-tech.ac.jp}}

\end{document}